\documentstyle[12pt]{article}
\begin{document}
\newtheorem{theorem}      {Th\'eor\`eme}[section]
\newtheorem{theorem*}     {theorem}
\newtheorem{proposition}  [theorem]{Proposition}
\newtheorem{definition}   [theorem]{Definition}
\newtheorem{e-lemme}        [theorem]{Lemma}
\newtheorem{cor}   [theorem]{Corollaire}
\newtheorem{resultat}     [theorem]{R\'esultat}
\newtheorem{eexercice}    [theorem]{Exercice}
\newtheorem{rrem}    [theorem]{Remarque}
\newtheorem{pprobleme}    [theorem]{Probl\`eme}
\newtheorem{eexemple}     [theorem]{Exemple}
\newcommand{\preuve}      {\paragraph{Preuve}}
\newenvironment{probleme} {\begin{pprobleme}\rm}{\end{pprobleme}}
\newenvironment{remarque} {\begin{rremarque}\rm}{\end{rremarque}}
\newenvironment{exercice} {\begin{eexercice}\rm}{\end{eexercice}}
\newenvironment{exemple}  {\begin{eexemple}\rm}{\end{eexemple}}
%
%
\newtheorem{e-theo}      [theorem]{Theorem}
\newtheorem{theo*}     [theorem]{Theorem}
\newtheorem{e-pro}  [theorem]{Proposition}
\newtheorem{e-def}   [theorem]{Definition}
\newtheorem{e-lem}        [theorem]{Lemma}
\newtheorem{e-cor}   [theorem]{Corollary}
\newtheorem{e-resultat}     [theorem]{Result}
\newtheorem{ex}    [theorem]{Exercise}
\newtheorem{e-rem}    [theorem]{Remark}
\newtheorem{prob}    [theorem]{Problem}
\newtheorem{example}     [theorem]{Example}
\newcommand{\proof}         {\paragraph{Proof~: }}
\newcommand{\hint}          {\paragraph{Hint}}
\newcommand{\heuristicproof}{\paragraph{heuristic proof}}
\newenvironment{e-probleme} {\begin{e-pprobleme}\rm}{\end{e-pprobleme}}
\newenvironment{e-remarque} {\begin{e-rremarque}\rm}{\end{e-rremarque}}
\newenvironment{e-exercice} {\begin{e-eexercice}\rm}{\end{e-eexercice}}
\newenvironment{e-exemple}  {\begin{e-eexemple}\rm}{\end{e-eexemple}}
\newcommand{\1}        {{\bf 1}}
\newcommand{\pp}       {{{\rm I\!\!\! P}}}
\newcommand{\qq}       {{{\rm I\!\!\! Q}}}
\newcommand{\B}        {{{\rm I\! B}}}
\newcommand{\cc}       {{{\rm I\!\!\! C}}}
\newcommand{\N}        {{{\rm I\! N}}}
\newcommand{\R}        {{{\rm I\! R}}}
\newcommand{\D}        {{{\rm I\! D}}}

\newcommand{\C}        {{\bf C}}        
\newcommand{\rank}{\hbox{rank}}
\newcommand{\CC}{{\cal C}}
\def \Re {{\rm Re\,}}
\def \Im {{ \rm Im\,}}
\def\Hom{{\rm Hom\,}}
\def\Lip{{\rm Lip}}
%
%
\newcommand{\dontforget}[1]
{{\mbox{}\\\noindent\rule{1cm}{2mm}\hfill don't forget : #1
\hfill\rule{1cm}{2mm}}\typeout{---------- don't forget : #1 ------------}}
\newcommand{\note}[2]
{ \noindent{\sf #1 \hfill \today}

\noindent\mbox{}\hrulefill\mbox{}
\begin{quote}\begin{quote}\sf #2\end{quote}\end{quote}
\noindent\mbox{}\hrulefill\mbox{}
\vspace{1cm}
}
\title{ Filling real hypersurfaces by pseudoholomorphic discs}
\author{ Alexandre Sukhov{*} and Alexander Tumanov{**}}
\date{}
\maketitle

{\small
{*}Universit\'e des Sciences et Technologies de Lille,
Laboratoire Paul Painlev\'e, U.F.R. de Math\'ematique,
59655 Villeneuve d'Ascq, Cedex, France,
sukhov@math.univ-lille1.fr

{**} University of Illinois, Department of Mathematics
1409 West Green Street, Urbana, IL 61801, USA, tumanov@math.uiuc.edu
}
\bigskip

Abstract. We study pseudoholomorphic discs with boundaries
attached to a real hypersurface $E$ in an almost
complex manifold. We give sufficient conditions for
filling a one sided neighborhood of $E$ by the discs.
\bigskip


\bigskip

\section{Introduction}

Since the fundamental work of Gromov \cite{Gr}
pseudoholomorphic curves have become an object of intensive
research because of their remarkable applications in symplectic
geometry and low dimensional topology; see \cite{AuLa, MS}
for numerous references.
Pseudoholomorphic discs with boundaries in a prescribed
real manifold play an important role. We call them the
{\it Bishop discs} after Bishop \cite{Bi}, who introduced
a method for constructing such discs.
We refer the reader to the expository papers \cite{MePo,Tr2,Tu2}
for the history and recent results on the Bishop
discs in the case of the standard complex structure.

The main result of the paper is the following

\begin{e-theo}
\label{Through}
Let $E$ be a $C^\infty$-smooth real hypersurface in an almost complex
manifold $(M,J)$. Assume that the set of all points where the Levi
form of $E$ vanishes identically has empty interior.
Also assume that $E$ contains no $J$-complex hypersurface
passing through a point $p\in E$. Then the Bishop discs of $E$
fill a one-sided neighborhood of $p$.
\end{e-theo}

This paper is a continuation of our previous work \cite{SuTu},
in which we prove a similar result for a manifold $M$ of complex
dimension 2; see \cite{SuTu} for additional discussions.
For the standard complex structure, a closely related result
was obtained by Tr\'epreau \cite{Tr}, who proved the forced
one-sided extendibility of CR-functions rather than filling
by the discs.
The second author \cite{Tu} proved an analogue for a generic
submanifold $E$ of arbitrary codimension in $\cc^n$.
We must stress that \cite{Tr, Tu} don't need any assumptions
on the Levi form.

In \cite{SuTu} we prove that if $M$ has complex dimension 2,
and the hypersurface $E\subset M$ does not contain any $J$-complex
curves, then the Bishop discs fill a one-sided neighborhood
of each point of $E$. We don't need the Levi non-flatness
assumption in \cite{SuTu}, because in complex dimension 2,
if $E$ is Levi-flat, then it is foliated by $J$-complex hypersurfaces,
that is, $J$-complex curves.
This is no longer true in complex dimension 3 or higher due to
an example by Ivashkovich and Rosay \cite{IvRo}.
In this example a real hypersurface $E$ in an almost complex
manifold of complex dimension 3 has identically vanishing Levi form
but $E$ contains no complex hypersurfaces. Thus, in general,
the vanishing of the Levi form of a real hypersurface does not imply
integrability of the distribution of holomorphic tangent spaces.
Nevertheless every Bishop disc of $E$ is contained in $E$,
which is always true for real hypersurfaces with identically
vanishing Levi form \cite{KrSu}. Hence the assumption on the Levi
form can not be dropped in Theorem \ref{Through}. It would be
interesting  to replace it by a weaker assumption.

Thus the non-integrability of an almost complex structure in complex
dimension $3$ or higher leads to a new phenomenon in the
behavior of Bishop discs. This is the main motivation of our work.

Our method is based on our previous work \cite{SuTu}.
For the convenience of the reader we use similar notations and
terminology.

The paper was written for a special volume of JGA
issued on the occasion of Gennadi Henkin's 65th birthday.
The second author expresses deep gratitude to Henkin for advising
him as a student, and for life-long support and inspiration.
Both authors wish Henkin good health, happiness, and
continuing success for many years to come.

\section{Preliminaries}
In this section we briefly recall some basic properties
of almost complex manifolds.
Throughout the paper, we assume by default that all relevant objects
(manifolds, structures, etc. ) are smooth of class $C^\infty$.

\subsection{Almost complex structures}
Recall that a real $2n$-dimensional manifold $M$ with
a fixed $(1,1)$-tensor field $J$ satisfying $J^2 = -Id$
is called {\it an almost complex manifold}; the tensor
field $J$ is called {\it an almost complex structure}.

A real submanifold $Y$ of an almost complex manifold
$(M,J)$ is called $J$-{\it complex} if for every
point $p \in Y$, the tangent space $T_pY$ is a $J$-complex
subspace of $T_pM$, that is, $J(p)T_pY=T_pY$.
A $J$-{\it complex hypersurface} is a $J$-complex
submanifold of real dimension $2n-2$.

We denote by $\D$  the unit disc in $\cc$ and by $J_{st}$
the standard complex structure in $\cc^n$.
 Let $(M,J)$ be an almost complex manifold and let $f$ be
 a smooth map from $\D$ into $M$. As usual we say that $f$ is
 {\it $J$-holomorphic}  if $df \circ J_{st} = J \circ df$.
 We call such a map $f$
a $J$-{\it holomorphic} disc or a {\it pseudoholomorphic} disc.

An  almost complex manifold
$(M,J)$ of complex dimension $n$ can be  locally viewed  as
the unit ball $\B$ in
$\cc^n$ equipped with a small almost complex
deformation of $J_{st}$. More precisely, given point $p \in
M$,  for every  $\delta_0 > 0$, and   every  $k
\geq 0$
 there is a neighborhood $U$ of $p$ and a
smooth coordinate chart  $z: U \longrightarrow \B$ such that
$z(p) = 0$, $dz(p) \circ J(p) \circ dz^{-1}(0) = J_{st}$,  and the
direct image $z_*(J) := dz \circ J \circ dz^{-1}$ satisfies
the inequality
$\vert\vert z_*(J) - J_{st}
\vert\vert_{C^k(\bar {\B})} \leq \delta_0$.

\begin{definition}
Consider  a $C^2$ function $u$   on $M$. Let $p$ be a point of $M$
and $v \in T_pM$ be a tangent vector. The Levi form of the
function $u$ at $p$ evaluated on $v$ is defined
by the equality $L^J(u)(p)(v):=-d(J^* du)(v,Jv)(p)$.
The Levi form of a real hypersurface
$E = \{ \rho = 0 \}\subset M$ at $p \in E$ is the conformal class of
the Levi form $L^J(\rho)(p)$ of the defining function $\rho$
on the $J$-holomorphic tangent space
$H^J_p E:= T_p E \cap J T_p E$.
\end{definition}

It is well known (see, e. g., \cite{IvRo}) that
$L^J(u)(p)(v) = \Delta(u \circ f)(0)$
where $f$ is an arbitrary $J$-holomorphic disc in $M$ such that
 $f(0) = p$ and $
df(0)(\partial_{\Re \zeta}) = v$ (here $\zeta$ is the
standard complex coordinate variable  in $\cc$). We also point
out that the Levi form is  invariant with respect to
$J$-biholomorphisms:
if $\Phi$ is a $(J,J')$-holomorphic diffeomorphism from $(M,J)$
into $(M',J')$,
then $L^J(u)(p)(v) =L^{J'}(u \circ \Phi^{-1})(\Phi(p))(d\Phi(p)(v))$.

We use only the {\it quadratic} Levi form defined above.
Similarly, one could define a {\it bilinear} version of the Levi form.
Then vanishing of the bilinear Levi form of a real hypersurface
$E$ is equivalent to Frobenius-integrability of
the $J$-holomorphic tangent bundle $H^J E$.
However, in the almost complex case, the bilinear Levi form
is not necessarily symmetric, so vanishing of the quadratic
Levi form does not imply that $E$ is foliated by $J$-complex
hypersurfaces.

We refer the reader to \cite{AuLa, Ch, MS, MePo} for further
details on almost complex structures, CR-manifolds, and
pseudoholomorphic discs.

\subsection{Normal form of an almost complex structure along a pseudoholomorphic disc}
Let $J$ be a smooth almost complex structure on a small enough neighborhood
$U$  of the origin in $\cc^n$
and $J(0) = J_{st}$.
Denote by $z = (z_1,\dots, z_n)$ the standard coordinates in $\cc^n $.
In the matrix computations throughout this paper we always view $z$
as a column. We will use the following notation:
$\overline\partial  := \frac{\partial}{\partial \overline \zeta},\;
\partial : = \frac{\partial }{\partial \zeta}$.

Then a map $z:\D \longrightarrow U$ is $J$-holomorphic if and only if
it satisfies the following  system of partial differential equations (the
Cauchy - Riemann equations):
\begin{eqnarray}
\label{basicequation}
\overline \partial z - A(z)\overline{\partial z} = 0,
\end{eqnarray}
where $A(z)$ is the complex $n\times n$ matrix  defined by
\begin{eqnarray}
\label{matrixA}
A(z)v = (J_{st} + J(z))^{-1}(J_{st} - J(z))(\overline v).
\end{eqnarray}
It is easy to see that the right-hand side
is $\cc$-linear in $v\in\cc^n$ with respect to the standard
structure $J_{st}$, hence $A(z)$ is well defined.
Since $J(0) = J_{st}$, we have $A(0) = 0$.  Furthermore,
by shrinking the neighborhood
of the origin where the initial disc $z$ is contained, we can
assume that given positive $k$ the $C^k$ norm of $A$ is
as small as we need. This simple argument is used
throughout the paper. We denote by $e_n$
the last vector of the standard basis of $\cc^n$ that is
$$
e_n = (0,\dots,1).
$$

\begin{e-lemme}
\label{normalization} Let $z_0$ be a $J$-holomorphic
disc  sufficiently close in the $C^k$ norm for some $k \geq 2$ to the disc
$\zeta \mapsto \zeta e_n$, $\zeta\in\D$.
Then there exists a local coordinate diffeomorphism in a neighborhood of
$z_0(\overline\D)$ such that in the new coordinates we have
$z_0(\zeta) = \zeta e_n$, $\zeta\in\D$.
Furthermore, in the new coordinates we have $A(\zeta e_n) = 0$,
$A_z(\zeta e_n) = 0$ for $\zeta \in \D$.
\end{e-lemme}

This statement is proved in \cite{SuTu} in the case of complex
dimension 2. The proof goes through in any dimension
with obvious modifications, so we omit it.

Let $\rho$ be a real function in a neighborhood of the origin
in the space  $\cc^n$ equipped with a smooth almost complex
structure $J$. Even if $J(0) = J_{st}$, the Levi for of $\rho$ with
respect to $J$ at the origin does not necessarily coincide with
the Levi form with respect to $J_{st}$. However, if the coordinates
are normalized according the previous lemma, the Levi forms of $\rho$
with respect to $J$ and $J_{st}$ are the same
(see, e. g., \cite{SuTu}).

\begin{e-lemme}
\label{Leviform}
Assume that $A(0) = A_z(0) = 0$. Then
the Levi form of $\rho$ at the origin with respect to the
structure $J$ coincides with the Levi form of $\rho$
at the origin with respect to the structure $J_{st}$.
\end{e-lemme}

\section{Bishop discs and the Bishop equation} Let $(M,J)$ be a
smooth almost complex
manifold of complex dimension
$n$ and $E$  a real generating submanifold of codimension $d$ in $M$. Recall that a $J$-holomorphic disc $f:\D \longrightarrow M$ continuous on
$\overline \D$ is called a {\it Bishop disc} if $f(b\D) \subset E$,
where $b\D$ denotes the boundary of $\D$.

If $p\in E$ and a Bishop disc $f$ additionally satisfies
the condition $f(1) = p$ then we say that $f$
{\it is attached to $E$ at} $p$.

Let
$\rho = (\rho_1,\dots,\rho_d)$ be a local defining function of $E$
in a neighborhood $U$ of a point $p$ where local coordinates $z$
are fixed.
Then a smooth map $z:\D \longrightarrow U$ continuous on
$\overline\D$ is a Bishop disc if and only if it satisfies
the system (\ref{basicequation}) of partial differential
equations with the boundary condition $(\rho \circ z)(\zeta) = 0$,
$\zeta \in b\D$. We call this boundary value problem
{\it the Bishop equation} by analogy with the case of $\cc^n$.
The neighborhood $U$ is supposed to be small enough so we call
such Bishop discs {\it small}. The Bishop equation can be easily
solved by the Implicit Function Theorem using functional properties
of the standard integral operators. This implies the existence
and local parametrization of Bishop discs attached to $E$.
The following existence result is obtained in \cite{SuTu}.

\begin{e-pro}
\label{theosolvingBishop}  Given positive non-integral $k$ the set ${\mathcal A}^J_p(E)$ of $J$-holomorphic Bishop discs of $E$ of class $C^k(\overline\D)$  satisfying $z(1) = p$
is a Banach manifold and its tangent space at the constant disc  $z  \equiv p$ is
canonically isomorphic to the  space $(C^k(\overline\D) \cap {\mathcal O}_1(\D))^{n-d}$ of vector functions valued in $\cc^{n-d}$, holomorphic (with respect to $J_{st}$) on $\D$, of class $C^k$ on $\overline \D$ and vanishing at the point $1$. The manifold  ${\mathcal A}^J_p(E)$ depends smoothly on a point $p$ and a deformation of the complex structure $J$.
\end{e-pro}

The following example is due to Ivashkovich and Rosay
(\cite{IvRo}, section 6).
\medskip

{\bf Example 1.}
Consider in $\cc^3$ with the standard complex coordinates
$z = (z_1, z_2, z_3)$ the almost complex structure $J$ such that
the equations (\ref{basicequation}) for a $J$-holomorphic disc
$z: \zeta \mapsto z(\zeta)$ have the form
\begin{eqnarray*}
& &\overline\partial z_j = 0,\; j = 1,2\\
& &\overline\partial z_3 =\overline z_1 \overline{\partial z_2}
\end{eqnarray*}
An equivalent definition of $J$ as a real matrix function is given
in \cite{IvRo}. It is shown in \cite{IvRo} that the real hypersurface
$E = \{ \Im z_3 = 0 \}$ contains no $J$-complex hypersurfaces.
On the other hand it is easy to see that the Levi form of $E$ with
respect to $J$ vanishes identically at every point and according
to the general result \cite{KrSu} all Bishop discs of $E$ are
contained in $E$. In fact, in this example the Bishop discs $z$
attached to $E$ at the origin can be described as follows.

Since $z$ is $J$-holomorphic, then $z_1$ and $z_2$
are holomorphic functions on the unit disc in the usual sense.
Then the equation for $z_3$ implies
$z_3=\phi+\bar\psi$, where $\phi$ is holomorphic and
$\psi(\zeta)=\int_1^\zeta z_1(\tau)z_2'(\tau)\,d\tau$.

Since $z$ is a Bishop disc for $E$ then $\Im z_3|_{b\D} = 0$.
Since $z_3$ is harmonic, then $\Im z_3$ vanishes on
$\D$ identically, that is, $z(\D)\subset E$.
Then
$z_3(\zeta)=2\Re\int_1^\zeta z_1(\tau)z_2'(\tau)\,d\tau$,
and $z_1$ and $z_2$ are arbitrary holomorphic functions in $\D$
with $z_1(1)=z_2(1)=0$.

This description implies that the evaluation map
$$
{\cal F}: {\cal A}^J_0(E)\ni z \mapsto
z(-1)=\left(z_1(-1), z_2(-1),
2\Re\int_1^{-1} z_1(\tau)z_2'(\tau)\,d\tau\right)\in E
$$
has the rank 5 except at the constant disc $z=0$.
It also implies that the range of ${\cal F}$ is all of $E$.
Thus, all Bishop discs attached to $E$ at the origin
lie in $E$ and fill all of it.

Example 1 shows that the notion of the {\it defect} of a
Bishop disc introduced in \cite{Tu} does not make sense
in the almost complex setting.
According to \cite{Tu}, for a small Bishop disc $z$ through
a fixed point $p=z(1)$ of a real generic manifold
$E\subset(\C^n, J_{st})$,
the spaces of infinitesimal perturbations of $z'(1)$
and those of $q=z(-1)$ have the same codimension in
the respective ambient spaces. This codimension is called the defect
of the disc $z$. In Example 1, for a non-constant disc $z$,
we found that the map $z\mapsto z(-1)$ has the rank 5
whereas the map $z\mapsto z'(1)$ has the rank at most 4 because
$z(\D)\subset E$. Hence the approach based on the defect
fails in this example.

Baouendi, Rothschild, and Tr\'epreau \cite{BaTrRo} interpret the
defect of a disc in terms of its lifts to the conormal
bundle of $E$ in the cotangent bundle $T^*\cc^n$.
Although an almost complex structure admits natural lifts to the
cotangent bundle of an almost complex manifold (see, e.g., \cite{YI}),
obviously, they don't give rise to the correct notions of the defect
of Bishop discs.
\bigskip

The following simple consequence of Proposition
\ref{theosolvingBishop} reduces the proof of the main result
to constructing a single transverse Bishop disc similarly
to \cite{SuTu, Tu}.

\begin{e-pro}
\label{filling}
Let $E$ be a real hypersurface in an almost complex manifold $(M,J)$. Suppose that there exists $z \in {\cal A}_p^J(E)$ such that the normal derivative vector
$dz(\partial_{\Re \zeta}\vert_1)$
is not tangent to $E$ at $p$ (that is the Bishop disc $z$ is attached to $E$ at $p$ transversally).
Then Bishop discs of $E$ fill a one-sided neighborhood of $p \in E$.
\end{e-pro}

In some special cases the condition of transversality can be easily
verified. For instance, if $E$ is a pseudoconvex hypersurface
or satisfies standard finite type conditions \cite{SuTu}.
However, this is a substantially more subtile task to prove
the existence of a transversal Bishop disc without these additional
assumptions. This is our main goal in the present paper.

Our key tool is the following statement.

\begin{e-theo}
\label{coordinatefree}
Let  $E$ be a generating submanifold in an almost complex manifold $(M,J)$ and let $f_0$ be a small enough embedded Bishop disc attached to $E$
at a point $p \in E$  and tangent to $E$ at $p$.
Suppose that every Bishop disc $f$ attached at $p$ and close enough to $f_0$
also is tangent to $E$ at $p$. Then the Levi form of every defining function of $E$ vanishes identically on the holomorphic tangent space of $E$
at every point of the boundary of the disc $f_0$.
\end{e-theo}

In the present paper we will use this proposition only in the case where $E$ is a hypersurface.

\begin{e-cor}
Let  $E$ be a real hypersurface in an almost complex manifold $(M,J)$ and let $f_0$ be a small enough embedded Bishop disc attached to $E$
at a point $p \in E$  and tangent to $E$ at $p$.
Suppose that every Bishop disc $f$ attached at $p$ and close enough to $f_0$
also is tangent to $E$ at $p$. Then the Levi form of $E$ vanishes identically
along the boundary of the disc $f_0$.
\end{e-cor}

In the hypothesis of Theorem \ref{coordinatefree} we will suppose that
$$
E = \{ \rho = (\rho_1,\dots,\rho_d) = 0 \}
$$
is a real generating $C^\infty$-smooth submanifold of codimension $d$
in an almost complex manifold $(M,J)$ of complex dimension n.
Since all our considerations are purely local, we work in local
coordinates. For technical reasons it is more convenient to consider
Bishop discs attached to $E$ at the point $e_n$ so we suppose that
$E$ contains it.  Our almost complex structure $J$ is viewed as
a real $(2n \times 2n)$-matrix function $J: z \mapsto J(z)$.
We also can assume  that $J(e_n) = J_{st}$.  We view a defining function $\rho$ of $E$ as a vector function valued in the space $\R^d$ of real columns. Recall also that a map $\zeta \mapsto z(\zeta)$ from
the unit disc $\D$ to $\cc^n$ is $J$-holomorphic if and only if it satisfies (\ref{basicequation}).
We stress that we consider here only maps valued in a small enough neighborhood of the point $e_n$.

Consider small embedded $J$-holomorphic Bishop discs attached to $E$ at the point $e_n$. For every  such a disc $z_0$ there exists a
local diffeomorphism  such that in the new coordinates we have
$z_0(\zeta) =  \zeta e_n$, $\zeta\in\D$
(see Lemma \ref{normalization}).
We denote again by $J$ the representation of our almost complex structure in the new coordinates,
$J(e_n) = J_{st}$.

We establish the following coordinate version of Theorem \ref{coordinatefree}.

\begin{e-theo}
\label{MainProp}
Suppose that for the Bishop disc
$z_0:\zeta \mapsto  \zeta e_n$, $\zeta\in\D$,
we have $A \circ z_0 = 0$ and $A_z \circ z_0 = 0$ (here  $A$ is the matrix
from the Cauchy - Riemann equations (\ref{basicequation})) and
\begin{eqnarray*}
& &\rho^k_{z_j}(e_n) = \delta^k_j, \; k,j = 1,...,d,\\
& &\rho^k_{z_j}(e_n) = 0, \; k = 1,...,d, \; j=d+1,..,n
\end{eqnarray*}
where $\delta^k_j$ denotes Kroneker's symbol.
Assume that the derivatives $A_{\overline z}$ and the second derivatives of the functions $\rho^k$ are small on the disc $z_0(\overline\D)$.
Suppose further that for every Bishop disc
$z: \zeta \mapsto  z(\zeta)$
with
$z(1) = z_0(1) =  e_n$
close enough to $z_0$ we have $$ \partial z_k(1) = 0, k = 1,...,d$$
Then $L^J(\rho^k)(p)(v)  = 0$, $k = 1,...,d$  for every point $p \in z_0(b\D)$ and every vector $v \in H^J_p(E)$.
\end{e-theo}

{\bf Remark 1.}
The assumptions of smallness of the derivatives of $A$ and
$\rho$ are automatically satisfied because the original
disc is small and the change of coordinates stretches it
to the unit disc.
The normalization condition  also can be achieved
by a suitable change of coordinates which does not affect previous
assumptions (Lemma \ref{normalization}). So Theorem \ref{coordinatefree} is a
  consequence of Theorem \ref{MainProp}.

{\bf Remark 2. }We point out that in the above coordinates the Levi form of
$\rho$ with respect to $J_{st}$ coincides with the Levi form with respect to
$J$ at any point of the set $z_0(b\D)$ (Lemma \ref{Leviform}).

We provide the proof of Theorem \ref{MainProp} in the next two sections.

\section{Infinitesimal analysis of the  Bishop equation}
In this section we linearize Bishop's equation and then solve
the obtained linear boundary problem for an elliptic system of
partial differential equations. This gives a convenient parametrization of the tangent space
$T_{z_0}{\cal A}^J_p(E)$, $p = e_n$. In particular this allows to reformulate
the moment condition of tangency of all Bishop discs attached at $p$ as a
condition of vanishing of a non-linear map $\Phi$ defined on the Banach manifold ${\cal
  A}^J_p(E)$. All the statements of this section are proved in \cite{SuTu} for the case of complex dimension 2 and
  the proofs remain valid in the general case with obvious changes. For the convenience of readers
  we outline here the main steps of this construction.

\subsection{Linearized Bishop equation}

Consider a Bishop disc
$z: \zeta \mapsto  z(\zeta)$, $\zeta\in\D$,
attached to $E$ at $e_n$ and close enough to $z_0$.
Then this disc satisfies the following boundary problem:
\begin{eqnarray*}
& &\overline\partial z = A \overline{\partial z},\\
& &\rho \circ z \vert_{b\D} = 0.
\end{eqnarray*}

Recall that  we only deal with small Bishop discs, so we can assume that for
every fixed non-integral $k > 0$
the norm $\parallel A \parallel_{C^k}$ is small enough.

Now we introduce solutions $\dot z$ of the corresponding linearized problem:
\begin{eqnarray*}
& &\overline\partial \dot z = \dot A(z) \overline{\partial z} + A
\overline{\partial \dot z},\\
& &\Re (\rho_z(z)\dot z)\vert_{b\D} = 0,
\end{eqnarray*}
where $\rho_z$ is the Jacobi matrix:
$$\rho_z = \left(
\begin{array}{ccc}
\rho^1_{z_1} &    \dots & \rho^1_{z_n}\\
\dots & \dots &    \dots\\
\rho^d_{z_1} &  \dots &  \rho^d_{z_n}
\end{array}
\right)$$
The symbol $\dot A$ is defined as follows.
Consider the map
$z \mapsto A \circ z$ defined on the space of smooth discs
$z :\D \longrightarrow \cc^n$.
Then $\dot A(z)$ denotes the value of the Frechet derivative
of this map at  $z$ and $\dot z$ (a tangent vector at $z$):
$\dot A = A_z\dot z + A_{\overline z} \overline{\dot z}$.
As usual we call the solutions $\dot z$ of this linearized system
the {\it infinitesimal variations  of the disc $z$}.

In order to obtain a convenient description of the space of
solutions of the linearized Bishop equation, following \cite{SuTu},
we make suitable changes of dependent variables. Namely, put
$I' = (I - A \overline A)^{-1}$ and
$\lambda = \rho_z I' + \rho_{\overline z} \overline{I'}\overline A$,
$$
\lambda = \left(
\begin{array}{ccc}
\lambda^1_{1} &    \dots & \lambda^1_{n}\\
\dots & \dots &    \dots\\
\lambda^d_{1} &  \dots &  \lambda^d_{n}
\end{array}
\right)
$$

Consider also the matrix $\Lambda$

$$\Lambda  = \left(
\begin{array}{ccccc}
\lambda^1_1 &  \dots & \lambda^1_{d+1} & \dots & \lambda^1_n\\
\dots & \dots &  \dots & \dots &  \dots\\
\lambda^d_{1} &  \dots & \lambda^d_{d+1} & \dots  & \lambda^d_{n}\\
0 & \dots & 1 &  \dots & 0\\
\dots & \dots & \dots & \dots &  \dots\\
0 &  \dots & 0 &\dots & 1
\end{array}
\right)$$
We assume that $\Lambda$ is smoothly extended
on the disc $\D$.
Introduce the change of variable
$$v = \Lambda(\dot z - A \overline{\dot z}),
\qquad
v = (v_1,\dots,v_n).$$
\begin{e-pro}
\label{change}
The new dependent variable $v$ has the following properties.
\begin{itemize}
\item[(i)] $v$ is a solution of the equation
$$\overline\partial v = B_1 v + B_2\overline v$$
where  $B_j = B_j(z)$ are certain smooth $(n \times n)$ - matrices.
\item[(ii)] $v$ satisfies the boundary condition
$$
\Re v = \left(
\begin{array}{cl}
0\\
u
\end{array}
\right), v(1) = 0
$$
where $u: b\D \longrightarrow \R^{n-d}$,
$u = (u_{d+1},\dots,u_{n})$, is an arbitrary smooth vector
function with $u(1) = 0$.
\item[(iii)] We have $$\partial v_k(1) = 0$$
for $k = 1,...,d$.
\item[(iv)]  $\dot B_1 - (\overline\partial \dot \Lambda) \Lambda^{-1}$  is a linear combination
of $\dot z$, $\overline{\dot z}$ with smooth
coefficients depending on $\zeta$,
which we write in the form
\begin{eqnarray}
\label{mod1}
\dot B_1 = (\overline\partial \dot \Lambda) \Lambda^{-1} \,\,\,
\hbox{mod} \, (\dot z, \overline{\dot z})
\end{eqnarray}
and similarly $\dot B_2$  is a linear combination
of $\dot z$, $\overline{\dot z}$
and $\overline\partial \overline{\dot z}$ with smooth
coefficients depending on $\zeta$:
\begin{eqnarray}
\label{mod2}
\dot B_2 = 0 \,\,\,
\hbox{ mod} \, (\dot z, \overline{\dot z}, \overline\partial \overline{\dot z}).
\end{eqnarray}
\end{itemize}
\end{e-pro}
The proof is given in \cite{SuTu}.

\subsection{Parametrization of the solution space of the  linearized Bishop equation}
In order to describe solutions of the linearized Bishop equation we need some
integral operators which we briefly recall here. The Cauchy-Green transform is defined by
\begin{eqnarray}
\label{CauchyGreen}
Tf(\zeta) = \frac{1}{2\pi i} \int\int_{\D} \frac{f(\tau)}{\tau -
  \zeta}d\tau \wedge d\overline\tau.
\end{eqnarray}
It is classical  (\cite{Ve}, p. 56, Theorem 1.32) that for every integral  $m
\geq 0$ and    $0 < \alpha < 1$ the linear map
$T:C^{m + \alpha}(\overline \D) \longrightarrow C^{m + \alpha +1}(\overline \D)$
is bounded.

We use the following notation for  the Cauchy integral

\begin{eqnarray}
\label{Cauchy}
Kf(\zeta) = \frac{1}{2 \pi i}\int_{b\D} \frac{f(\tau)d\tau}{\tau - \zeta}.
\end{eqnarray}

We also set
\begin{eqnarray*}
K_1u: = Ku - (Ku)(1)
\end{eqnarray*}

The boundary problem stated in Proposition \ref{change} (i),
(ii) was solved in \cite{SuTu}.  Its solution  has the form
\begin{eqnarray}
\label{parametrization}
v = v_0 + R_1 v_0 + R_2\overline v_0,
\end{eqnarray}
Here the resolvent operators $R_j$ are
$\cc$-linear bounded operators
$C^{m + \alpha} \longrightarrow C^{m + \alpha + 1}$
for any integral $m \geq 0$ and   $0 < \alpha < 1$  and
$$v_0 =  \left(
\begin{array}{cl}
0\\
\varphi
\end{array}
\right),$$
where $\varphi$ is a holomorphic vector function given by
$$\varphi = 2K_1 u.$$

One can view $R_j$ as $(n \times n)$-matrices with operator entries.
Furthermore, in \cite{SuTu} an expansion of $R_1$ to a sum of operator
series is obtained and studied. Every term of this series is
a composition of several integral operators, which are suitable
modification of the Cauchy-Green operator $T$, and the operators
of left multiplication by the matrices $B_j$, $\overline B_j$.

Now we can  interpret the  condition
\begin{eqnarray}
\label{*}
 \partial v_1(1) = \dots = \partial v_d(1) = 0
\end{eqnarray}
(see (iii) in Proposition \ref{change}) similarly to \cite{SuTu}.
We have
$$v_j = \sum_{k=d + 1}^n R_1^{jk}\varphi_k + \sum_{k=d+1}^n R_2^{jk}\overline\varphi_k$$
for $j = 1,...,d$. Here $R_1^{jk}$ and $R_2^{jk}$ denote the (jk) matrix entries of the matrix operators $R_1$ and $R_2$ respectively.
Thus the condition (\ref{*}) means that
\begin{eqnarray*}
\sum_{k=d+1}^n \partial R_1^{jk}\varphi_k(1) + \sum_{k=d+1}^n \partial R_2^{jk}\overline\varphi_k(1) = 0
\end{eqnarray*}

Since the first term is $\cc$-linear and the second is anti-linear, we get $\sum_{k=d+1}^n \partial R_1^{jk}\varphi_k(1)  = 0$
for all $\varphi$ with $\varphi(1) = 0$. Hence,

\begin{eqnarray}
\label{moment}
 \partial R_1^{jk}\varphi_k(1)= 0, \; j = 1,...,d, \; k = d+1,...,n
\end{eqnarray}
for all holomorphic functions $\varphi_k$ with $\varphi_k(1) = 0$.

Now we rewrite the  condition (\ref{moment}) in a more
convenient form. The method of \cite{SuTu} is based on the smoothing
properties of the operators involved to the expansion of $R_1$.
Differentiation of these operators in (\ref{moment})  leads to certain singular integral
operators whose  kernels have singularities at the point $1$. These operators loose the smoothing
property at this point.
In \cite{SuTu} we overcome this technical difficulty  by introducing suitable
function spaces.  The idea is to consider spaces of functions which are smooth
on $\D$ except the point $1$. Then we also need to
adapt the notion of a bounded linear operator to this class of  spaces.

Given integral $m \geq 0$ and $0 < \alpha < 1$
the spaces $C^{m + \alpha}_1(\overline\D)$ and $C^{m+\alpha}_1(b\D)$
are defined as spaces of functions which are of class
$C^{m + \alpha}$ on  $\overline\D \backslash \{ 1 \}$.
More precisely
$f \in C_1^{m+\alpha}(\overline\D)$ if $f \in L^\infty(\D)$
and for every $\varepsilon > 0$ we have
$f\vert_{\overline\D \backslash \B(1,\varepsilon)}
\in C^{m +\alpha}(\overline\D \backslash \B(1,\varepsilon))$,
where $\B(1,\varepsilon)$ denotes the disc of radius
$\varepsilon$ centered at $1$.
The space $C_1^{m+\alpha}(b\D)$ is defined similarly.
As in \cite{SuTu}
we say that $P$ is {\it a bounded linear operator}
$C^{m+\alpha}_1(\overline\D) \longrightarrow
C_1^{k+\beta}(\overline\D)$ for integral $m,k \geq 0$ and  $0 < \alpha, \beta < 1$,
if for every $\varepsilon > 0$ there exists a constant
$C = C(\varepsilon)>0$ such that for every
$f \in C^{m+\alpha}_1$ we have
\begin{eqnarray*}
\parallel Pf \parallel_{C^{k+\beta}(\D \backslash
\B(1,2\varepsilon))}
+\parallel P f \parallel_{L^\infty(\D)} \leq
C (\parallel f \parallel_{C^{m+\alpha}(\D \backslash
\B(1,\varepsilon))}
+\parallel f \parallel_{L^\infty(\D)}).
\end{eqnarray*}
Similarly are defined  bounded operators
$C^{m+\alpha}(\overline\D) \longrightarrow
C_1^{k+\beta}(\overline\D)$ and
when we have $b\D$ in place of $\D$.
\begin{e-def}
\label{equivalence}
We write  $P_1 \sim P_2$ for two operators
$P_1$ and $P_2$ if $P_1 - P_2$ is a bounded operator
$C^{1+\alpha}(\overline\D) \cap {\cal O}(\D) \longrightarrow
C_1^{2+\alpha}(b\D)$. Here ${\cal O}(\D)$ denotes the space of holomorphic (in
the usual sense) functions on $\D$.
\end{e-def}

The proof of the following three propositions is given in \cite{SuTu}. They
form  crucial  steps of our reduction.

First, we reduce the condition (\ref{moment}) to a vanishing of  certain
non-linear map on the space of Bishop discs.

\begin{e-pro}
\label{Derivation1}
 The condition (\ref{moment}) is equivalent to the fact
that for any
$j = 1,...,d$ and $k = d+1,...,n$
we have
\begin{eqnarray}
\label{functional}
\Phi^{jk}(z):= T\left ( B_1 +  \Omega\right )^{jk} \vert_{b\D} = 0
\end{eqnarray}
where a map $\Omega: z  \mapsto  \Omega(z)$ associates to $z$
an $(n \times n)$-matrix with operator entries and the index $(jk)$
refers to the matrix entry.
This  property holds for every Bishop disc $z$ close enough to $z_0$
with $z(1) =  e_n$.
\end{e-pro}
 A precise description and analysis of the map
$\Omega$ is contained in \cite{SuTu}. Here we just point out
that in the notation of \cite{SuTu} the this map  admits an
expansion $\Omega = \overline B_2 \mu(\tau) + \Sigma' + \Sigma''$,
where $\Sigma'$ and $\Sigma''$ are the sums of certain operators
series. Every term of these series is a composition of certain
singular integral operators and operators of the left multiplication
by the matrices $B_j$, $\overline B_j$. These series converge
in the space $C^{\1-\beta}(\overline\D)$ for every
$0 < \beta < 1$.

 The map $\Phi = (\Phi^{k}) $  is defined  in
a neighborhood of the disc $z_0$ on  a Banach manifold
${\cal A}_p^J(E)$ of $J$-holomorphic Bishop discs attached to $E$
at the point $p = e_n$. It follows from  (\ref{functional})
that  the map $\Phi$ vanishes identically. Hence  its Fr\'echet
derivative $\dot \Phi$ at $z_0$  also does. This leads to the second step (the
linearization):

\begin{e-pro}
\label{Differentiation2}
 For the Fr\'echet derivative of $\Phi$ we have
\begin{eqnarray}
\label{Frechet}
\dot\Phi = (\dot \Phi^{jk}) = \left ( T\left [ \dot B_1 + \dot \Omega \right ]^{jk} \vert_{b\D} \right ) =0.
\end{eqnarray}
where as above $j = 1,...,d$, $k = d+1,...,n$.
\end{e-pro}

 In other words  $\dot\Phi(\dot z) = 0$,
for every $ \dot z \in T_{z_0}{\cal A}_p^J(E)$.

  In view of the parametrization of the tangent
space $T_{z_0}{\cal A}_p^J(E)$ given by (\ref{parametrization})
we can view   $\dot \Phi$ as an $\R$-linear operator applied to
a vector-function
$\varphi \in ({\cal O}(\D) \cap C^{1+\alpha}(\overline \D))^{n-d}$,
$0 < \alpha < 1$.

\begin{e-pro}
\label{Derivation3}
The condition (\ref{Frechet}) implies that
\begin{eqnarray}
\label{Step1}
T(I_0\dot B_1)^{jk}\vert_{b\D}\sim 0,
\end{eqnarray}
$j = 1,...,d$, $k = d+1,...,n$.
Here $I_0 = I + a_1 + a_2 \mu$,
$\mu(\tau) = \left ( \frac{\tau - 1}
{\overline \tau - 1} \right )^2$, and each matrix function
$a_j$ is of class $C^{m + \alpha}_1(\overline\D)$ for all integral
$m \geq 0$, $0 < \alpha < 1$.
\end{e-pro}

 The matrices $a_j$ here are written
in \cite{SuTu} in the terms of the expansion of $\Omega$ mentioned
after Proposition \ref{Derivation1}. The matrices $a_j$ are small, so $I_0$ is close to $I$.
In particular the matrix $I_0$ is invertible. The equivalence
$\sim$ in (\ref{Step1}) is understood in the sense of
Definition \ref{equivalence}.

\section{Infinitesimal analysis of the operator  $\Phi$}
Now we are able to prove Theorem \ref{MainProp}.
The idea is to remove the terms contained in $\Omega$
in (\ref{Frechet}). The equivalence $\sim$
of operators is always understood in the sense of
Definition \ref{equivalence}.

Consider $[T,b]:= T \circ b - b \circ T$  the commutator
of the operator $T$ and the operator of multiplication by a matrix $b$.
The next technical statement also is due to \cite{SuTu}.
\begin{e-lemme}
\label{commutator}
Let $b\in C^{m+\alpha}_1(\overline\D)$ for any  all integral $m \geq 0$, $0 < \alpha < 1$.
Then $[T,b]$ defines a bounded operator
$C^s(\overline \D) \longrightarrow C^{s+2}_1(\overline \D)$
for any non-integral $s > 0$.
\end{e-lemme}

In other words $T \circ b \sim b \circ T$, which we use below.
The first step is the following
\begin{e-pro}
We have
\begin{eqnarray}
\label{Frechet3}
(T\dot B_1)^{jk}\vert_{b\D}\sim 0.
\end{eqnarray}
for $j=1,...,d$, $k = d+1,...,n$.
\end{e-pro}
\proof Using the commutator  argument and Lemma \ref{commutator} we obtain from (\ref{Step1}) that
$$(I_0 T \dot B_1)^{jk}\vert_{b\D}\sim 0$$
for $j = 1,...,d$ and $k = d+1,...,n$. Then
$$
0\sim(I_0 T \dot B_1)^{jk}\vert_{b\D}
= \sum_{s=1}^n I_0^{js} T\dot B_1^{sk}.
$$
Since
$$
\dot \Lambda =   \left(
\begin{array}{ccccc}
\dot \lambda^1_1 &  \dots & \dot \lambda^1_d & \dots & \dot \lambda^1_n\\
\dots & \dots &  \dots & \dots &  \dots\\
\dot \lambda^d_{1} &  \dots & \dot \lambda^d_d & \dots  & \dot \lambda^d_{n}\\
0 & \dots & 0 &  \dots & 0\\
\dots & \dots & \dots & \dots &  \dots\\
0 &  \dots & 0 &\dots & 0
\end{array}
\right)$$
and by (\ref{mod1}) we have
\begin{eqnarray*}
\dot B_1 = (\overline\partial \dot \Lambda)\Lambda^{-1} \,\,\,
\hbox{mod} \,
(\dot z, \overline{\dot z}),
\end{eqnarray*}
we obtain that
$$T\dot B_1^{sk}\vert_{b\D}\sim 0$$ for $s = d+1,...,n-1$ and
every $k$.
Therefore
$$\sum_{s=1}^d I_0^{js}(T\dot B_1)^{sk}\vert_{b\D}\sim 0$$
for $j = 1,...,d$, $k = d+1,...,n$.
Recall that the  matrix $I_0$ is close to $I$ (see the end of Section 4 after
Proposition \ref{Derivation3}), so the
matrix $(I_0^{js})_{j,s = 1,...,d}$ is invertible.
This implies (\ref{Frechet3}).
\bigskip

On the disc $z_0$ we have
$$
\lambda \circ z_0=  \left(
\begin{array}{ccc}
\rho^1_{z_1} &   \dots & \rho^1_{z_n}\\
\dots & \dots &   \dots\\
\rho^d_{z_1} &  \dots &  \rho^d_{z_n}\\
\end{array}
\right)
$$

$$
\Lambda \circ z_0=  \left(
\begin{array}{ccccc}
\rho^1_{z_1} &  \dots &  \rho^1_{d+1} & \dots & \rho^1_{z_n}\\
\dots & \dots &  \dots & \dots &  \dots\\
\rho^d_{z_1} &  \dots & \rho^d_{d+1} & \dots & \rho^d_{z_n}\\
0 & \dots &  1 & \dots &  0\\
\dots & \dots &  \dots & \dots & \dots\\
0 &  \dots  & 0 & \dots & 1
\end{array}
\right)
$$
Using the obvious block structure of the above matrix, we have
$$
\Lambda \circ z_0=  \left(
\begin{array}{cc}
X & Y\\
0 &  I
\end{array}
\right)
\qquad
(\Lambda \circ z_0)^{-1} =
\left(
\begin{array}{cc}
X^{-1} &  -X^{-1}Y\\
0 &  I
\end{array}
\right),
$$
where $X = (\rho^j_{z_k})_{j,k = 1,...,d}$.
Now using the condition (\ref{Frechet3})
and  Lemma \ref{commutator}, we conclude that
\begin{eqnarray}
\label{Tdbar}
T(\overline\partial\dot\lambda)\vert_{b\D}
\left(
\begin{array}{cl}
-X^{-1}Y\\
I
\end{array}
\right)
\sim 0,
\end{eqnarray}

Finally we need to express $\dot \lambda$ in terms of $\varphi$.
We use the notation $\lambda^k$ for the row
$\lambda^k = (\lambda^k_1,..., \lambda^k_n)$ of the matrix $\lambda$.
The condition  $A_z\circ z_0=0$ implies
$$
\dot \lambda^k = (\rho^k_z)\dot{} + \rho^k_{\overline z}
\dot{\overline A} = (a^k\dot z + b^k\overline{\dot z})^t,
$$
where
$a^k = \rho^k_{zz}+\rho^k_{\overline z}\overline A_z$ and
$b^k = \rho^k_{z\overline z}$ and the upper index $t$ denotes
the matrix transposition.

We have $\dot z = \Lambda^{-1} v_0 + (S)$ where $(S)$ denotes
smoother terms. Then

\begin{eqnarray*}
\dot z
\sim \Lambda^{-1} \left(
  \begin{array}{cl}0\\
  \varphi
  \end{array}
\right)=
\left(
\begin{array}{cc}
X^{-1} &  -X^{-1}Y\\
0 &  I
\end{array}
\right)
\left(
\begin{array}{c}
0\\
\varphi
\end{array}
\right)
\end{eqnarray*}
By the Cauchy-Green formula, $T\overline\partial\varphi = 0$ and $T\overline\partial\overline\varphi = \overline\varphi$. Again by Lemma \ref{commutator} we obtain
$$
T\overline\partial\dot\lambda^k\vert_{b\D}
\sim
(a^kT\overline\partial\dot z
+b^kT\overline\partial\overline{\dot z})^t
\sim
(b^k\overline{\dot z})^t
\sim
(\overline{\dot z_1},\dots, \overline{\dot z_n})
 \left(
\begin{array}{ccc}
\rho^k_{z_1\overline{z}_1}&\dots &\rho^k_{z_n \overline{z}_1}\\
\dots & \dots & \dots\\
\rho^k_{z_1\overline{ z}_n} & \dots &\rho^k_{z_n \overline{z}_n}
\end{array}
\right)
$$
Then the conditions (\ref{Tdbar}) can be rewritten in the form
$$
\overline\varphi^t L(\rho^k)  \sim 0, k = 1,...,d
$$
where
\begin{eqnarray*}
 L(\rho^k) = \left(
\begin{array}{cc}
-(\overline{X}^{-1}\overline Y)^t  & I
\end{array}
\right)
 \left(
\begin{array}{ccc}
\rho^k_{z_1\overline{z}_1}& \dots &\rho^k_{z_n \overline{z}_1}\\
\dots & \dots  & \dots\\
\rho^k_{z_1\overline{z}_n} &  \dots &\rho^k_{z_n \overline{z}_n}
\end{array}
\right)
\left(
\begin{array}{c}
 -X^{-1}Y\\
I
\end{array}
\right)
\end{eqnarray*}
Thus the product $ \overline \varphi^t L(\rho^k)$ is of class
$C^{2+\alpha}_1(b\D)$ for any vector function
$\varphi \in (C^{1+\alpha}(\overline\D) \cap {\cal O}(\D))^{n-d}$.
Hence the matrix $L(\rho^k)$ vanishes identically on the boundary
of the disc $z_0$. But this is exactly the matrix of the restriction
of the Levi form of the function $\rho^k$
on the holomorphic tangent space (with respect to $J_{st}$)
at a point of $z_0(b\D)$.

Recall that the conditions $A \circ z_0 = 0$ and
$A_z \circ z_0 = 0$ imply that the Levi form of $\rho^k$
with respect to the structure $J$ coincides with the Levi
form of $\rho^k$ with respect to the structure
$J_{st}$ at every point of the boundary of the disc $z_0$ (Lemma \ref{Leviform}).
Thus the Levi form with respect to the structure $J$ of every defining function of $E$
 vanishes on the boundary of the disc $z_0$, as desired.

This proves Theorem \ref{MainProp}.

\section{The case of degenerate rank} We return to the case where
$E = \{ \rho = 0 \}$ is a hypersurface. In this section we study
the situation where the boundaries of the
Bishop discs attached to $E$ through a fixed point do not cover
an open set in $E$. The method of \cite{SuTu} can be easily
adapted to this case.

\begin{e-pro}
\label{MainProp2}
Suppose that the boundaries of $J$-holomorphic discs
$\zeta\mapsto z(\zeta)$ with
$z(1) =  e_n$
attached to $E$ and close to the disc
$z_0:\zeta \mapsto  \zeta e_n$, $\zeta\in\D$,
do not cover an open set in $E$. Then for every $\zeta_0 \in b\D$,
$\zeta_0 \neq 1$ there
exists a $J$-complex hypersurface near   the point $\zeta_0 e_n$
completely contained in $E$.
\end{e-pro}

\proof Fix a point $\zeta_0 \neq 1$ in $b\D$. For every Bishop disc $z(\zeta)$ close enough to $z_0(\zeta)$
we define the evaluation map  ${\cal F}_{\zeta_0}: z \mapsto z(\zeta_0)$.
Its Fr\'echet derivative $\dot {\cal F}_{\zeta_0}$  at $z_0$ is given by
$\dot {\cal F}_{\zeta_0}: \dot z \mapsto \dot z(\zeta_0)$
where $\dot z$ is an infinitesimal perturbation of $z_0$.

By assumption the boundaries of discs do not cover an open
subset of $E$. Therefore
$\rank \dot {\cal F}_{\zeta_0} \leq 2n-2$ for all $\zeta_0$.
First we show  that $\rank\,\dot {\cal F}_{\zeta_0} \geq 2n-2$.

\begin{e-lemme}
\label{lemma8}
For every vector
$(q_2,\dots,q_n) \in \cc^{n-1}$ there exists
$\dot z =(\dot z_1,\dot z_2,\dots,\dot z_n)$
with $(\dot z_2(\zeta_0),\dots,\dot z_n(\zeta_0))
 =( q_2,\dots, q_n)$.
\end{e-lemme}
\proof Recall again that $A = 0$ and $A_z = 0$ on the disc $z_0$ and
\begin{eqnarray*}
& &\dot z = \Lambda^{-1}v,\\
& &v = \left(
\begin{array}{cl}
0\\
\varphi
\end{array}
\right) + R_1 \left(
\begin{array}{cl}
0\\
\varphi
\end{array}
\right) + R_2 \left(
\begin{array}{cl}
0\\
\overline\varphi
\end{array}
\right),
\end{eqnarray*}
where $\varphi = (\varphi_2,\dots,\varphi_n)$ is an arbitrary
holomorphic vector function valued in $\cc^{n-1}$ satisfying
$\varphi(1) = 0$. Then for $j = 2,...,n$
\begin{eqnarray*}
\dot z_j = \varphi_j + \sum_{k=2}^n R_1^{jk}\varphi_k
+ \sum_{k=2}^n R_2^{jk}\overline \varphi_k.
\end{eqnarray*}

 Plugging $\zeta = \zeta_0$,
we get for $j=2,...,n$
\begin{eqnarray}
\label{**}
\dot z_j(\zeta_0) = \varphi_j(\zeta_0)+
\sum_{k=2}^n\int\int_{\D}a_{jk}^1(\zeta)\varphi_k(\zeta)
d\zeta \wedge d\overline\zeta
+\sum_{k=2}^n\int\int_{\D}a_{jk}^2(\zeta)\overline{\varphi_k(\zeta)}
d\zeta \wedge d\overline\zeta.
\end{eqnarray}
where $a_{jk}^i$  are integrable functions in $\D$ (this integral
representation follows from the analysis of the expansion of $R_j$
into operator series in \cite{SuTu}).  If the rank of the map
$\varphi \mapsto (\dot z_2(\zeta_0),\dots,\dot z_n(\zeta_0))$
is smaller than or equal to $2n-3$, then there exists a vector
$c = (c_2,...,c_n) \in \cc^{n-1} \backslash \{ 0 \}$ such that
for every $\varphi$ we have
$Re(\sum_{j=2}^{n}c_j \dot z_j(\zeta_0))= 0$.
Then for some $b_{j}^1, b_{j}^2 \in L^1(\D)$ we have
\begin{eqnarray*}
2\Re \left ( \sum_{j=2}^n c_j\varphi_j(\zeta_0) \right )
+ \sum_{j = 2}^n \int\int_{\D} b_{j}^1(\zeta)\varphi_j(\zeta)
d\zeta \wedge d\overline\zeta
+\sum_{j=2}^n\int\int_{\D} b_{j}^2(\zeta)\overline{\varphi_j(\zeta)}
d\zeta \wedge d\overline\zeta = 0.
\end{eqnarray*}
Splitting into linear and anti-linear parts we get
\begin{eqnarray*}
\sum_j \left ( c_j\varphi_j(\zeta_0)
+ \int\int_{\D}b_{j}^1(\zeta)\varphi_j(\zeta)
d\zeta \wedge d\overline\zeta \right ) = 0.
\end{eqnarray*}
Similarly to \cite{SuTu} this implies that $c_j = 0$ for every
$j = 2,...,n$. Indeed, take for instance $\varphi_j = 0$,
$j=3,...,n$ and $\varphi_2 = \psi^n$, where $\psi$ has a peak
at $\zeta_0$, that is, $\psi(0) = 0$, $\psi(\zeta_0) = 1$ and
$\vert \psi(\zeta) \vert < 1$ for
$\zeta \in \overline \D \backslash \{ \zeta_0 \}$. Then passing to
the limit as $n \longrightarrow \infty$ we obtain that $c_2 = 0$.
We proceed similarly for other $c_j$. The obtained contradiction
proves the lemma.
\bigskip

Thus the rank of the map $\varphi \mapsto \dot z(\zeta_0)$
is equal to $2n-2$. Hence there exist vectors $d^1,d^2 \in \cc^{n-1}$
such that
\begin{eqnarray}
\label{***}
\dot z_1(\zeta_0) = \sum_{j= 2}^{n}d^1_j \dot z_j(\zeta_0)
+ \sum_{j=2}^{n} d^2_j\overline{\dot z_j(\zeta_0)}
\end{eqnarray}
 for all $\varphi$.
The equality
$\dot z = \Lambda^{-1}v$ implies
\begin{eqnarray}
\label{****}
\dot z_1 = -\rho_{z_1}^{-1}\rho_{z_2}\varphi_2 -...
-\rho_{z_1}^{-1}\rho_{z_n} \varphi_n + \sum_k P_{1k}\varphi_k
+ \sum_k P_{2k}\overline\varphi_k,
\end{eqnarray}
where $P_{1k}$ and $P_{2k}$ are integral operators.
Expressing $\dot z(\zeta_0)$ in terms of $\varphi$ by (\ref{*}),
(\ref{****}), we get from (\ref{***})
\begin{eqnarray*}
& &-\rho_{z_1}^{-1}\rho_{z_2}\varphi_2(\zeta_0) -...
-\rho_{z_1}^{-1}\rho_{z_n} \varphi_n(\zeta_0)
= \sum d^1_j\varphi_j(\zeta_0)
+ \sum_j d^2_j\overline{\varphi_j(\zeta_0)}\\
& &+
\sum_k\int\int_{\D}b_{1k}(\zeta)\varphi_k(\zeta)
d\zeta \wedge d\overline\zeta
+\sum_k\int\int_{\D}b_{2k}(\zeta)\overline{\varphi(\zeta)}d\zeta
\wedge d\overline\zeta
\end{eqnarray*}
for some $b_{1k},b_{2k} \in L^1(\D)$.
As in lemma \ref{lemma8} we obtain
$$
d^1_j = -\rho_{z_1}^{-1}\rho_{z_j}\vert_{\zeta = \zeta_0}, j=2,...,n
$$
$$
d_j^2 = 0, j = 2,...,n
$$

Since $\zeta_0 \in b\D$ is arbitrary, we
have $$\sum_{j=1}^{n} \rho_{z_j}\dot z_j  = 0$$ on $b\D$. This
precisely means that  $\dot z(\zeta) \in H^J_{z(\zeta)}E$ (the holomorphic
tangent space) for any $\vert \zeta \vert = 1$.

By the hypothesis of proposition this is true for every disc close to
$z_0$.
By the Rank Theorem, the image of the evaluation map ${\cal F}_{\zeta_0}$ is a
$J$-complex hypersurface contained in $E$. This completes the proof of the proposition.

\subsection{Proof of the main result}

It is convenient to begin with the proof of a weaker version of our main result.

\begin{e-pro}
\label{MainTheorem}
Let $E$ be a real hypersurface in an almost complex manifold $(M,J)$.
Assume that the set of points where the Levi form of $E$ vanishes
identically has the empty interior and $E$ contains no $J$-complex
hypersurfaces. Then  Bishop discs of $E$ fill a one-sided neighborhood
of every point of $E$.
\end{e-pro}

\proof
If $E$ admits a transversal Bishop  disc attached at $p$, then the statement follows
by Proposition \ref{filling}. Suppose that there are no transversal Bishop discs.
 Then by Theorem \ref{MainProp}
the Levi form of $E$ (with respect to $J$)  vanishes on the boundary of every
Bishop disc attached to $E$ at $p$. If these discs fill an open
subset $\Omega$ of $E$, the Levi form of $E$ vanishes on $\Omega$ identically which contradicts to the assumption of theorem.
 Finally, if the boundaries of Bishop
discs do not cover an open piece of $E$, then Proposition \ref{MainProp2}
implies the existence of $J$-complex hypersurfaces in $E$.
Thus, $E$ necessarily admits a transversal Bishop disc which proves Proposition \ref{MainTheorem}.
\bigskip

{\bf Proof of Theorem  \ref{Through}:}
The proof goes along the lines of that of Proposition \ref{MainTheorem}.
We only describe the necessary adjustments.
We assume that there is no transverse Bishop disc attached to
$E$ at $p$ and bring this to a contradiction.

To obtain a complex hypersurface passing through $p=e_n$
in Proposition \ref{MainProp2} we use a disc $z_0$ with
$z_0(\zeta_0)=z_0(1)$, $\zeta_0=-1$. To obtain such a disc
we take an arbitrary embedded Bishop disc $f$ and put
$z_0(\zeta)=f(\zeta^2)$.

Therefore in proving most auxiliary results, we have to
deal with Bishop discs $z$ close to the disc
$z_0:\zeta\mapsto \zeta^2 e_n$.
Such discs generally cannot be straightened by
a diffeomorphism in the proof of Lemma \ref{normalization}
on the normalization of the matrix $A$ along a disc,
so we need a suitable version of the lemma.
In the new version we require that
the conditions $A=0$, $A_z=0$ be satisfied on
$z(\overline\D_0)$, where $\D_0\subset\D$ is a fixed
subdomain so that $b\D_0$ contains an open arc
$\gamma\subset b\D$, $-1\in\gamma$, and the map
$\zeta\mapsto\zeta^2$ is a diffeomorphism on $\D_0$.
We also require that $A(e_n)=0$.
The needed version easily follows from existing
Lemma \ref{normalization}.

With the above normalization of $A$, the proof of all
the results will go through with obvious adjustments.
Some relations will hold only on $\D_0$ or $\gamma$.
For instance, the formulas (\ref{mod1}) and (\ref{mod2})
will hold only on $\D_0$. In the equivalence relation in
Definition \ref{equivalence} we will change the target space
to $C_1^{2+\alpha}(\gamma)$.
In the proof of Proposition \ref{MainProp}  we
obtain that the Levi form of $E$ vanishes on
$z_0(\gamma)$.

Theorem \ref{Through} is proved.

\end{document}